\magnification=\magstephalf

\vbadness=10000
\parindent=15pt 
\hsize = 15truecm
\vsize = 21truecm 
\baselineskip=13pt
\parskip 1pt plus 1pt
\widowpenalty=1000 \clubpenalty=1000  

\nopagenumbers
\headline={\ifodd\pageno\rightheadline\else\leftheadline\fi}
\def\rightheadline{\pfont\hfil\number\pageno}
\def\leftheadline {\pfont\number\pageno\hfil}


\def\Ai{{{\rm Ai}}}
\def\maplenl{\par\ \ \ }

\font\tfont=cmss12 at 14pt
\font\afont=cmss12
\font\bfont=cmss9 
\font\sffont=cmss10

\font\cfont=cmssi9 

\font\reffont=cmcsc10
\font\titfont = cmbx10 scaled \magstep3
\font\chapfont = cmbx10 scaled \magstep2
\font\secfont=cmbx10 scaled \magstep1

\def\subsecfont{\bf}
\def\subsubsecfont{\sl}

\font\pfont=cmr10 at 12pt




\newcount\frnum                   
\newcount\exnum                   
\newcount\chapnum                 
\newcount\secnum                  
\newcount\subsecnum
\newcount\subsubsecnum
\newcount\fignum
\newcount\thnum
\newcount\remnum
\newcount\vbnum
\newcount\defnum
\newcount\lemnum

\def\title#1{ 
\vskip 2cm
\cl{\titfont#1}
\vskip 1cm
 \chapnum=0 \secnum=0 \exnum=0  \frnum=0 \fignum=0
 \thnum=0  \remnum=0 \vbnum=0 \defnum=0
     \message{#1} 
	\noindent\rm}

\def\chap#1{ 
\advance\chapnum by 1
\leftline{\chapfont \the\chapnum. #1}
 \secnum=0 \exnum=0  \frnum=0 \fignum=0
 \thnum=0  \lemnum=0 \remnum=0 \vbnum=0 \defnum=0
 \equationnumber=0
     \message{#1}  \def\chaptitle{#1}    
	\noindent\rm}

  \def\sect#1{\advance\secnum by 1\subsecnum=0 
               \vskip 20pt
               \leftline{\secfont \the\secnum. #1}
               \message{#1}\def\sectitle{#1}\nobreak\noindent}
                 
\def\subsect#1{\advance\subsecnum by 1\subsubsecnum=0
               \vskip 15pt
         \leftline{\subsecfont 
         \the\secnum.\the\subsecnum. #1}
               \message{#1}\nobreak \noindent}
\def\subsubsect#1{\advance\subsubsecnum by 1
               \vskip 10pt
         \leftline{\subsubsecfont 
         \the\secnum.\the\subsecnum.\the\subsubsecnum. #1}
               \message{#1}\nobreak \noindent}

\def\ex{\advance \exnum by 1
         \ifnum\exnum > 1 \vskip\baselineskip \fi
         \noindent {\bf \the\exnum.}\quad
         }

\def\fr{\advance \frnum by 1
         \ifnum\frnum > 1 \vskip\baselineskip \fi
         \noindent {\bf \the\frnum.}\quad
         }

\def\fig{\advance \fignum by 1
{{\bf Figure \the\fignum}\sq
}
         }

\def\remark{\advance \remnum by 1
         \vskip\baselineskip
         \noindent {\bf Remark \the\remnum.}\quad
         }

\def\example{\advance \vbnum by 1
         \vskip\baselineskip
         \noindent {\bf Example \the\vbnum.}\quad
         }

\def\defi{\advance \defnum by 1
         \vskip\baselineskip
         \noindent {\bf Definition \the\defnum.}\quad
         }

\def\theorem{\advance \thnum by 1
         \vskip\baselineskip
         \noindent {\bf Theorem \the\thnum.}\quad\sl
         }

\def\lemma{\advance \lemnum by 1
         \vskip\baselineskip
         \noindent {\bf Lemma  \the\lemnum.}\quad\sl
         }

\def\proof{\vskip\baselineskip\rm\noindent{\bf Proof.}\quad}
\def\eoproof{{\unskip\nobreak\hfil\penalty50
	\hskip2em\hbox{}\nobreak\hfil\vrule height4pt width5.5pt depth2pt
	\parfillskip=0pt\medbreak}}


\def\makebox#1#2#3
{\hbox{\vrule
       \vbox to  #1{\hrule \vss
                   \hbox to #2{\hss#3\hss}\vss
                   \hrule}\vrule}}

\def\ref#1{{\reffont#1}}

\def\sq{\quad}

\def\cl{\centerline}
\def\el{\eqalign}

\def\C{{\cal C}}
\def\iy{\infty}

\def\br#1#2{{{#1}\over{#2}}}

\def\brg#1#2{{\textstyle{{#1}\over\smash{#2}}}}

\def\sn{\sum_{n=0}^\iy}

\def\RR{{{\rm I}\!{\rm R}}}
\def\CC{\hbox{\rlap{$\,\,
  $\hbox{\vrule height6.2pt width.35pt depth-0.1pt}}$\rm C$}}

\def\phase{{\rm ph}}

\def\ell{elliptic }

\def\leaderfill{\leaders\hbox to 2em{\hss.\hss}\hfill}

\def\dddot#1{%
\setbox0=\hbox{$#1$}
\vbox to\baselineskip{\vss\hbox{%
 ${\mkern1mu.\mkern-2mu.\mkern-2mu.}
                \phantom{\char'177}$}
      \kern-\ht0
      \copy0}}

\def\insil#1{}

%
\newcount\equationnumber	
\def\eqnum{\relax
	\global\advance\equationnumber by 1
	\equationnumberformat{\the\equationnumber}%
	}%
\def\eqname#1{\relax
	\count255=\equationnumber
	\assignnumber{EN#1}\equationnumber
	\global\equationnumber=\count255
	\global\advance\equationnumber by 1
	\ifnum\csname EN#1\endcsname=\equationnumber
	\else
		\message{The equation number for ``#1'' is incorrect!}%
	\fi
	\equationnumberformat{\csname EN#1\endcsname}%
	}%
\def\equationnumberformat#1{\eqno(\equationnumbertype{#1})}%
\def\equationnumbertype#1{\number#1\relax}%
\def\referenceequation#1{\relax
	\assignnumber{EN#1}\equationnumber
	\equationnumbertype{\csname EN#1\endcsname}%
	}%
\def\forwardreferenceequation#1#2{\relax
	\global\advance\equationnumber by #2
	\assignnumber{EN#1}\equationnumber
	\global\advance\equationnumber by -1
	\global\advance\equationnumber by -#2
	\referenceequation{#1}%
	}%
%
\def\assignnumber#1#2{\relax
	\ifnum0<0\csname#1\endcsname
	\else
		\global\advance#2 by 1
		\expandafter\expandafter\expandafter
			\xdef\csname#1\endcsname{\the#2}%
	\fi
	}%

\def\en{\eqname}
\def\req#1{({\referenceequation{#1}})}

\def\bxf#1{
$$
\vbox{\hrule\hbox{\vrule\kern3pt
\vbox{\kern3pt\hbox{\kern3pt$\displaystyle{#1}\kern3pt$}\kern3pt}
\vrule}\hrule}
$$
}
%
%
\def\bxfn#1#2{
$$
\el{
\vbox{\hrule\hbox{\vrule\kern3pt
\vbox{\kern3pt\hbox{\kern3pt$\displaystyle{#1}\kern3pt$}\kern3pt}
\vrule}\hrule}}
\eqno(#2)
$$
}
\newcount\cochapnum
\newcount\cosectnum 
\newcount\cosubsnum 

\font\cochapfont = cmbx10 scaled \magstep1
\def\cosectfont{\rm}
\def\cosubsectfont{\sl}
\def\cochap#1{\advance\cochapnum by 1\cosectnum=0
               \leftline{\ \cochapfont %
               \the\cochapnum.\quad   #1}
               }
\def\cosection#1{\advance\cosectnum by 1\cosubsnum=0
               \leftline{\cosectfont 
               \quad\quad\the\cochapnum.\the\cosectnum.\    #1}
               }
\def\cosubsection#1{\advance\cosubsnum by 1
               \leftline{\cosubsectfont
               \quad\quad\quad\the\cochapnum.\the\cosectnum\the\cosubsnu.\  
               #1}
               }
\cochapnum=0  

\def\({\left(}
\def\){\right)}

\def\[{\left[}
\def\]{\right]}

\def\sk{\sum_{k=0}^\infty\,}

\bigskip

\tfont
\cl {Symbolic Evaluation of Coefficients in}
\bigskip
\cl {Airy-type Asymptotic Expansions}
\bigskip\bigskip \bigskip
\afont
\cl{Raimundas Vidunas$^{1,} {}^2$  and Nico M. Temme$^2$}
\medskip
\cfont
\cl{$^1$
Korteweg de Vries Instituut voor Wiskunde, 
University of Amsterdam, }
\cl{
Plantage Muidegracht 24, 
1018 TV Amsterdam, 
The Netherlands}
\medskip
\cl{$^2$
CWI,
P.O. Box 94079,
1090 GB Amsterdam,
The Netherlands}
\medskip
\cfont
\cl{e-mail: \tt vidunas@wins.uva.nl, vidunas@cwi.nl, nicot@cwi.nl}\medskip
\centerline               {}

\bfont
\parindent=25pt
{\pfont ABSTRACT}\par\noindent
{\narrower\noindent
Computer algebra algorithms are developed for evaluating the coefficients
in Airy-type asymptotic expansions that are obtained from 
integrals with a large parameter. The coefficients are defined from recursive 
schemes obtained from integration by parts. An application is given for 
the Weber parabolic cylinder function.

\vskip 0.8cm \noindent
\cfont
1991 Mathematics Subject Classification:
\bfont
41A60, 33C10, 33C15, 33F05, 65D20.
\par\noindent
\cfont
Keywords \& Phrases:
\bfont
Airy-type asymptotic expansions,
Maple algorithms, 
parabolic cylinder functions.
\par\noindent
\cfont
Note:
\bfont
Supported by NWO, project number 613-06-565.
\bigskip}
\rm
\parindent=15pt
\sect{Introduction}
\def\sk{\sum_{k=0}^\infty\,}   
\pageno=1

When constructing uniform asymptotic expansions of solutions of
differential equations or of functions defined by integrals, usually a
difficulty arises when the coefficients of the expansion are
constructed.  As shown in \ref{Olver} (1974) for the Airy-type
expansions of Bessel functions, recursion relations for the coefficients
can be obtained for the case that the expansion is obtained by using a
linear second order differential equation.

In many publications this method has been used, for example in 
\ref{Olver} (1959) and \ref{Dunster} (1989), and for expansions
involving Bessel functions or parabolic cylinder functions similar
results are available. Having such a recursion relation for the
coefficients does not always give the possibility to obtain analytic
expressions of a number of coefficients,  because the recursion involves
integrals of previous coefficients together with a function that is not
easy to handle. Sometimes the coefficients can be explicitly expressed
in terms of coefficients of simpler expansions because different types
of expansions may be valid in overlapping domains. See for the Bessel
functions the relations in \ref{Olver} (1974), page 425, Exercise 10.3 or
\ref{Abramowitz \& Stegun} (1964), page 368, formula 9.3.40.

For special functions usually the same type of uniform expansions can be
obtained by using integral representations of the functions. Sometimes,
in a particular problem, the integral is the only tool available for
constructing uniform expansions. By using transformations of variables in
the integrals, these representations can be transformed into standard
forms for which an integration by parts procedure can be used to obtain
expansions in terms of, for example, Airy functions. 

Although it is usually not possible to derive recursion relations for the
coefficients obtained in this way, in all cases for special functions
known so far, it is possible to construct a number of coefficients, and
only because of the complexity of the problem, which implies limitations
with respect to available computer memory when doing symbolic
computations,  there is an upper bound for this number. An advantage of
the differential equation approach is the possibility to construct
realistic and sharp error bounds for the remainders in the expansions;
similar bounds cannot be obtained in the approach based on integral
representations.

In this paper we use integral representations and give Maple algorithms
for constructing the coefficients in uniform asymptotic expansions
involving Airy functions. First we describe how to obtain the
coefficients for a general case. For an application we obtain the
coefficients for the case of a special function called parabolic
cylinder function. Straightforward computations are often complicated
by appearance of algebraic roots in the output or intermediate expressions.
These algebraic roots can be avoided by replacing some
parameters with algeraic expressions in suitable new variables. In the
example of parabolic cylinder function we avoid computations with algebraic
roots by using variable $u$ (instead of $t$), and at the end we simplify
the output by introducing variable $\xi$. In the last section we give the
Maple code used for this example.

\sect{Airy-type asymptotic expansions}%
We consider integrals of the form
$$F_\eta(z)=\br1{2\pi i}\int_{{\cal
C}}e^{z(\br13t^3-\eta t)}f(t)\,dt,\en{a1}$$  
where the contour is starts at infinity with $\phase t=-\pi/3$ and returns to 
infinity
with $\phase t=\pi/3$. We assume that the function $f(t)$ is analytic in the
neighbourhood of the contour. The parameter $z$ is large positive number, and 
$\eta$
is also assumed to be real. Extension to complex values of the parameters is
possible, but this will not be discussed in this paper.

In the case $f(t)=1$ we obtain the Airy function (\ref{Temme} (1996), 
page 101)
$$\br1{2\pi i}\int_{{\cal
C}}e^{z(\br13t^3-\eta t)}\,dt=z^{-\br13}\Ai\left(\eta
z^{\br23}\right).\en{a2}$$
For more general functions $f$ the asymptotic expansion of $F_\eta(z)$ can be 
given
in terms of this Airy function. The asymptotic feature of this type of
integral  is
that the  phase function
$\phi(t)=\br13t^3-\eta t$ has two saddle points at
$\pm\sqrt{{\eta}}$ that coalesce when $\eta\to0$, and it is not possible to 
describe the
asymptotic behaviour of $F_\eta(z)$ in terms of simple functions when $\eta$ is
small. When the parameter
$\eta$ is positive and bounded away from 0, one can perform a saddle point 
analysis
on \req{a1} and use a conformal mapping
$\phi(t)-\phi(\sqrt{{\eta}})=\br12u^2$ with the condition $u(\sqrt{{\eta}})=0$. 
We obtain
$$F_\eta(z)=\br1{2\pi i}\int_{-i\iy}^{i\iy} e^{\br12zu^2}g(u)\,du,$$
where $g(u)=f(t)\,dt/du$, with  $dt/du=u/(t^2-\eta)$, which is regular at the
positive saddle point, but not at the negative saddle point. 
It follows that, when $\eta$ becomes small, a singularity due to $dt/du$ in the
$u-$plane approaches the origin, and an expansion of $dt/du$ at $u=0$
will have coefficients that become infinite as $\eta\to0$. Hence, by using
the standard saddle point method we obtain an expansion that is not
uniformly valid as $\eta\to0$. \par 

A modification of the saddle point
method is possible by taking into account both saddle points. We give an
integration by parts procedure that is a variant of Bleistein's method 
introduced in
\ref{Bleistein} (1966) (for a different class of integrals), 
and that gives the requested uniform expansion.

We assume that $f$ is an analytic function in a certain domain $G$ and write
$$f(t)=\alpha_0  +\beta_0t+(t^2-\eta)g(t),\en{a3}$$ where
$$
\alpha_0=\br12\left[f(\sqrt{{\eta}})+f(-\sqrt{{\eta}})\right],\sq \beta_0= 
\br1{2\sqrt{{\eta}}}\left[f(\sqrt{{\eta}})-f(-\sqrt{{-\eta}})\right].\en{a4}$$
Clearly $\alpha_0\to f(0), \beta_0\to f'(0)$ as $\eta\to0$ and the following
Cauchy integral representations hold
$$\el{
f(t) &=\br1{2\pi i}\int_\C\br{f(s)}{s-t}\,ds,\sq
\ g(t) =\br1{2\pi i}\int_\C\br{   f(s)}{(s-t)(s^2-\eta)}\,ds,\cr
\alpha_0    &=\br1{2\pi i}\int_\C\br{s f(s)}{s^2-\eta}\,ds,\sq
\beta_0    =\br1{2\pi i}\int_\C\br{   f(s)}{s^2-\eta}\,ds,\cr}
$$
where the contours of integration $\C$ encircle the points $t$ 
and/or $\pm\sqrt{\eta}$. 
Upon substituting \req{a3} in \req{a1}, we obtain
$$F_\eta(z)=z^{-\br13}\Ai\left(\eta z^{\br23}\right)\alpha_0-
z^{-\br23}\Ai'\left(\eta z^{\br23}\right)\beta_0+
\br1{2\pi i}\int_{{\cal
C}}e^{z(\br13t^3-\eta t)}(t^2-\eta)g(t)\,dt.$$
An integration by parts gives
$$F_\eta(z)=z^{-\br13}\Ai\left(\eta z^{\br23}\right)\alpha_0-
z^{-\br23}\Ai'\left(\eta z^{\br23}\right)\beta_0
-\br1{2\pi i}\int_{{\cal
C}}e^{z(\br13t^3-\eta t)}f_1(t)\,dt,$$
where $f_1(t)=g'(t)$. Repeating this procedure we obtain the compound
expansion
$$F_\eta(z)\sim z^{-\br13}\Ai\left(\eta
z^{\br23}\right)\sn(-1)^n\br{\alpha_n}{z^n}- z^{-\br23}\Ai'\left(\eta
z^{\br23}\right)\sn(-1)^n \br{\beta_n}{z^n},\en{a5}$$
where the coefficients $\alpha_n, \beta_n$ are defined as in \req{a4} with
the function $f$ replaced with $f_n$, which in turn is defined by
the scheme
$$f_{n+1}(t)=g'_n(t),\sq f_n(t)=\alpha_n 
+\beta_n t+(t^2-\eta)g_n(t),\en{a6}$$
with  $n=0,1,2,\ldots$ and  $f_0(t)=f(t)$. The expansion in \req{a5} is
valid for large values of $z$ and holds uniformly with respect to $\eta$
in a neighbourhood of the origin. A more precise formulation can be
given, but more information can be found in the literature; see \ref{Olver} 
(1974) and
\ref{Wong} (1989).

The functions $f_n(t)$  defined in \req{a6} can be represented
in the form of Cauchy-type integrals. We have the following theorem.

\theorem 
Let the rational functions $R_n(s,t,\eta)$ be defined by
$$R_0(s,t,\eta)=\br1{s-t}, \sq R_{n+1}(s,t,\eta)= \br{-1}{s^2-\eta}\,\br{d}{ds}
R_n(s,t,\eta),\sq n=0,1,2,\ldots,\en{a7}$$
where $s, t, \eta\in\CC, s\ne  t, s^2\ne\eta$. Let $f_n(t)$ be defined by the
recursive scheme \req{a6},
where $f_0$ is a given analytic function in a domain $G$. Then we have
$$f_n(t)=\br1{2\pi i}\int_\C R_n(s,t,\eta) f_0(s)\,ds,$$
where $\C$ is a simple closed contour in $G$ that encircles the points $t$ and
$\pm\sqrt{\eta}$.
\proof The proof starts with
$$f_n(t)=\br1{2\pi i}\int_\C R_0(s,t,\eta) f_n(s)\,ds,$$
and in this representation the recursion relation \req{a6} for the functions 
$f_n$ is
used. More details can be found in
\ref{Olde Daalhuis \& Temme} (1994).
\eoproof

For the coefficients $\alpha_n, \beta_n$ we have a similar representation:
$$\alpha_n=\br1{2\pi i}\int_\C A_n(s,\eta) f_0(s)\,ds,
\sq
\beta_n=\br1{2\pi i}\int_\C B_n(s,\eta) f_0(s)\,ds,\en{a8}$$
where  $\C$ is a simple closed contour in $G$ that encircles the points 
$\pm\sqrt{\eta}$ and where $A_n(s,t)$ and $B_n(s,t)$ follow the same recursion 
\req{a7} as the rational functions
$R_n(s,t,\eta)$, with initial values
$$A_0(s,\eta)=\br {s}{s^2-\eta},\sq B_0(s,\eta)=\br{1}{s^2-\eta}.$$

We see that the coefficients $\alpha_n, \beta_n$ that play a role in the 
expansion
\req{a5} are well defined from an analytical point of view. However, from a
computational  point of view it may be quite difficult to evaluate the 
coefficients.  For a simple rational function like $f_0(t)=1/(t+1)$ the
computations are rather straightforward, and we can even use residue calculus
to evaluate the integrals  in \req{a8}:
$$\alpha_n= -A_n(-1,\eta),\sq \beta_n=-B_n(-1,\eta).$$ 
The first few values are in this case
$$
\el{
\alpha_0=-\br{1}{\eta-1},   \sq\sq 
& \beta_0=\br{1}{\eta-1}, \cr
\alpha_1=\br{\eta+1}{(\eta-1)^3},  \sq\sq 
& \beta_1= -\br{2}{(\eta-1)^3}, \cr
\alpha_2=-4 \br{2\eta+1}{(\eta-1)^5},  \sq\sq 
& \beta_2= 2\br{\eta+5}{(\eta-1)^5}, \cr
\alpha_3=4\br{2\eta^2+21\eta+7}{(\eta-1)^7},  \sq\sq 
& \beta_3= -40\br{\eta+2}{(\eta-1)^7}, \cr
\alpha_4=-280\br{\eta^2+4\eta+1}{(\eta-1)^9}, \sq\sq 
& \beta_4= 40\br{\eta^2+19\eta+22}{(\eta-1)^9}, \cr
\alpha_5=280\br{\eta^3+29\eta^2+65\eta+13}{(\eta-1)^{11}}, \sq\sq 
& \beta_5= -1120\br{2\eta^2+14\eta+11}{(\eta-1)^{11}}. \cr 
}
$$

For a more complicated or general function $f_0(t)$ even computer algebra 
manipulations give complicated expressions which are very difficult to
evaluate.
In the next section we develop an algorithm for computing the coefficients 
$\alpha_n, \beta_n$ when the values of the derivatives of 
$f_0(t)$ at $t=\pm\sqrt{\eta}$ are available.

\sect{How to compute the coefficients $\alpha_n, \beta_n$}%
We explain how the coefficients $\alpha_n, \beta_n$ of \req{a5} can be
computed. To avoid the square roots in the formulas we replace $\eta$ with 
$b^2$,
and we write
\req{a6} in the form
$$f_0(t)=f(t),\sq f_{n+1}(t)=g'_n(t),\sq f_n(t)=\alpha_n 
+\beta_nt+(t^2-b^2)g_n(t),$$
for $n=0,1,2,\ldots$. We assume that the function $f$ is analytic in a domain
$G$, that the series expansions  used in this section are
convergent in $G$, and that the points $\pm b$ are inside $G$. 
Furthermore, we assume the coefficients $p_k^{(1)}, p_k^{(2)}$ of 
the expansions
$$f(t)=\sk p_k^{(1)}(t-b)^k,\sq f(-t)=\sk p_k^{(2)}(t-b)^k\en{a9}$$
are available. 

\theorem {\bf Algorithm.\hfil\break}
Let coefficients $f_k^e, f_k^o$ be defined by
$$f_k^e=\br12\[p_k^{(1)}+p_k^{(2)}\],\sq 
f_k^o=\br12\[p_k^{(1)}-p_k^{(2)}\],\sq k=0,1,2\ldots$$
and coefficents $f_k^{o,e}$ 
by the recursion
$$b f_k^{o,e}=f_k^o-f_{k-1}^{o,e},\sq k\ge0,$$
with $f_{-1}^{o,e}=0$.
Next, define coefficients $\gamma_k, \delta_k$ by
$$\gamma_0=f_0^e,\sq \delta_0=f_0^{o,e},$$
and for $k\ge1$:
$$\el{
\gamma_k&=\sum_{j=1}^k \br{(-1)^{k-j}\,j\,(2k-j-1)!}
{(2b)^{2k-j}\,k!\,(k-j)!}\,f_j^e,\cr
\delta_k&=\sum_{j=1}^k \br{(-1)^{k-j}\,j\,(2k-j-1)!}
{(2b)^{2k-j}\,k!\,(k-j)!}\,f_j^{o,e}.\cr}
\en{a10}$$
Finally, let for $n\ge0$  coefficients $\gamma_k^{(n)}, \delta_k^{(n)}$ be
defined by the recursion
$$\el{
\gamma_k^{(n+1)}&=(2k+1)\delta_{k+1}^{(n)}+2b^2(k+1)\delta_{k+2}^{(n)},\cr
\delta_k^{(n+1)}&=2(k+1)\gamma_{k+2}^{(n)},\sq k=0,1,2,\ldots \ ,\cr}
\en{a11}
$$
with $\gamma_k^{(0)}=\gamma_k, \delta_k^{(0)}=\delta_k$.
Then the coefficients $\alpha_n, \beta_n$ of expansion \req{a5} are given by
$$\alpha_n=\gamma_0^{(n)},\sq \beta_n=\delta_0^{(n)},\sq n \ge 0.$$

\proof
The coefficients $f_k^e, f_k^o$ occur in the expansions
$$f_e(t)=\sk f_k^e(t-b)^k\sq f_o(t)=\sk f_k^o(t-b)^k,$$
where $f_e(t), f_o(t)$ are the even and odd parts of $f$:
$$f_e(t)=\br12[f(t)+f(-t)],\sq f_o=\br12[f(t)-f(-t)],$$
and the coefficients $f_k^{o,e}$ occur in the expansion
$$\br1t f_o(t)=\sk f_k^{o,e}(t-b)^k,$$
The coefficients $\gamma_k,\delta_k$ occur in the expansion
$$f(t)= \sk \gamma_k (t^2-b^2)^k+t\sk \delta_k (t^2-b^2)^k.\en{a12}$$
Observe that
$$f_e(t)=\sk \gamma_k(t^2-b^2)^k,\sq f_o(t)=t\,\sk \delta_k(t^2-b^2)^k,$$
and we will verify the first relation of \req{a10}. We write
$$\gamma_k=\br1{2\pi i}\int f_e\(\sqrt{z+b^2}\)\,\br{dz}{z^{k+1}},$$
where the contour is a small circle around the origin. Also,
$$\gamma_k=\br1{2\pi i}\int f_e(t)\,\br{2t\,dt}{(t+b)^{k+1}(t-b)^{k+1}},$$
where the contour is a small circle around $t=b$.

Substitute the expansion $f_e(t)=\sum_{j=0}^\infty f_j^e(t-b)^j$. Then,
$$\gamma_k=\sum_{j=0}^k f_j^e\ \br1{2\pi i}
\int \br{2t\,dt}{(t+b)^{k+1}(t-b)^{k+1-j}}.\en{a13}$$
Expand
$$ \br{2t}{(t+b)^{k+1}}=\sum_{m=0}^\infty q_m(t-b)^m.\en{a14}$$
We find, by using (\ref{Temme} (1996), p. 108)
$$(1-z)^{-a}=\sum_{n=0}^\infty\br{(a)_m}{m!}\,z^m=
\sum_{n=0}^\infty{-a\choose m}(-z)^m,$$
$$q_m=(-1)^m\br{(k-m)\,(k+m-1)!}{(2b)^{k+m}\,m!\,k!}.$$
When we use \req{a14} in \req{a13}, we only need $q_m$ with $m=k-j$. 
This gives the first result of \req{a10}.
The proof for $\delta_k$ is the same, because $(1/t)f_o(t)$ is
again even.  

The coefficients $\gamma_k^{(n)}, \delta_k^{(n)}$ are used in
$$f_n(t)= \sk \gamma_k^{(n)} (t^2-b^2)^k+t\sk \delta_k^{(n)}
(t^2-b^2)^k,$$ 
and the recursions in \req{a11} are easily verified, as is the final 
relation
$$\alpha_n=\gamma_0^{(n)},\sq \beta_n=\gamma_0^{(n)},\sq n \ge 0.$$
\eoproof

The first few values are of the coefficients $\gamma_k, \delta_k$ of the 
expansion in \req{a12} are
$$
\el{
\gamma_0= f_0^e, \sq\sq  
&\delta_0=  \br1b f_0^o,\cr
\gamma_1=\br1{2b}f_1^e , \sq\sq   
&\delta_1= \br1{2b^3}(bf_1^o-f_0^o),\cr
\gamma_2=\br1{2b^3} (2bf_2^e-f_1^e),  \sq\sq  
&\delta_2= \br1{8b^5}(2b^2f_2^o-3bf_1^o+3f_0^o),\cr
}
$$
and we observe, as in \req{a10}, negative powers of $b$. From a computational 
point of view, this may cause numerical instabilities, 
because the coefficients are analytic functions of
$b$ at $b=0$. For example, taking $f(t)=1/(t+1)$ again, we obtain 
$$\gamma_k=\br1{(1-b^2)^{k+1}},\sq \delta_k=-\br1{(1-b^2)^{k+1}},
\sq k=0,1,2,\ldots,$$
which follows from
$$\br1{t+1}=\br{1-t}{1-t^2}=\br{1-t}{(1-b^2)-(t^2-b^2)}=
\sk\br{(t^2-b^2)^k}{(1-b^2)^{k+1}}-
t\,\sk\br{(t^2-b^2)^k}{(1-b^2)^{k+1}}.$$

 From the representations in, for example, \req{a10}, 
we conclude that if we apply the 
algorithm for computing the coefficients $\alpha_n, \beta_n$ of expansion 
\req{a5},
starting with numerical values of the  coefficients $p_k^{(1)}, p_k^{(2)}$ of
\req{a9}, we may
encounter numerical instabilities when $b$ is small. For this reason, it is 
important to use exact values of $p_k^{(1)}, p_k^{(2)}$, and computer algebra is
of great help here. In the next section we consider a non-trivial case in which
obtaining the exact values of the coefficients $p_k^{(1)}, p_k^{(2)}$ of 
\req{a9} 
also needs special care.

\remark
In order to compute the coefficients $\alpha_n, \beta_n$ for $n=0,1,\ldots, N$ 
from the relation
$$\alpha_n=\gamma_0^{(n)},\sq \beta_n=\delta_0^{(n)},\sq n \ge 0$$
and the recursion in \req{a11}, we need the starting values for this recursion
$\gamma_k, \delta_k$ for $k=0,1,\ldots, 2N$. Hence, as follows from \req{a10}, 
we also need $p_k^{(1)}, p_k^{(2)}, k=0,1,\dots, 2N$ in the expansions in 
\req{a9}.

\sect{Application to parabolic cylinder functions}%
Weber parabolic cylinder functions are solutions of the differential
equation
$$\br{d^2y}{dx^2}-\left(\brg14x^2+a\right)y=0.\en{a15}$$
Airy-type expansions for the solutions of this equation can be found in 
\ref{Olver}
(1959), and are obtained by using the differential equation. In
this section we show how to obtain an integral representation like
\req{a1}, and how to apply the algorithm of the previous section for
deriving an Airy-type asymptotic expansion. 

A standard solution of \req{a15} is the integral
(see formula (19.5.4), page 688 in \ref{A\&S})
$$U(a,x)=\br{e^{\br14x^2}}{i\sqrt{2\pi}}
\int_\C e^{-xs+\br12s^2}s^{-a-\br12}\,ds,\en{a16}$$
where the contour $\C$ is a vertical line in the complex plane with $\Re s>0$. 

We consider large negative values of $a$, and use Olver's notation
$$a=-\brg12\mu^2,\sq x=\mu t\sqrt{2}.$$
Changing the variable of integration by writing $s\to\mu s/\sqrt{2}$,
we obtain
$$
U\left(-\brg12\mu^2,\mu t\sqrt{2}\right)=\br{e^{\br12\mu^2t^2}}{i\sqrt{2\pi}}
\left(\br\mu 2\right)^{\br12\mu^2+\br12}\int_\C e^{z\phi(s)} s^{-1/2}\,ds,$$
where
$$\phi(s)=\brg12s^2-2st+\ln s,\sq z=\brg12\mu^2.$$
The saddle points are obtained from the equation $\phi'(s)=0$, that is, from
$$\br{s^2-2st+1}{s}=0,$$
which gives two solutions
$$s_\pm=t\pm\sqrt{t^2-1}.$$
The saddle points coalesce when $t\to \pm 1$. Observe that in the new
variables the
differential equation \req{a15} transforms into
$$\br{d^2y}{dt^2}-\mu^4\left(t^2-1\right)y=0,$$
which has turning points at $t=\pm 1$.

A transformation into the standard form \req{a1} can be obtained by writing
$$\phi(s)=\brg13 w^3-\eta w + A,\en{a17}$$
where $\eta$ and $A$ have to be determined and do not depend on $t$.
A transformation into the cubic polynomial is first considered in
\ref{Chester} et
al. (1957).  For further details on the theory of this method we refer to
\ref{Olver} (1974), \ref{Wong} (1989), and \ref{Olde D. \& T} (1994).

The parameters $\eta$ and $A$ are obtained by assuming that the saddle points
$s_\pm$ in the $s-$variable should correspond with the saddle points
$w_\pm=\pm\sqrt{\eta}$ in the $w-$variable. We write
$$
t=\cosh \theta,\sq {\rm which\ gives}\sq
s_\pm=e^{\pm\theta},\en{a18}$$ 
assuming for the time being that $\theta\ge0$. 
We obtain the equations
$$\el{
\brg12e^{+2\theta}-2e^{+\theta}\cosh\theta+\theta&=-\brg23\eta^{3/2}+A,\cr
\brg12e^{-2\theta}-2e^{-\theta}\cosh\theta-\theta&=+\brg23\eta^{3/2}+A,\cr
}
$$
from which we derive
$$\brg43\eta^{3/2}=\sinh2\theta-2\theta,\sq 
A=-\brg12-\cosh^2\theta=-\brg12-t^2.\en{a19}$$
By using these values of $\eta$ and $A$ the $w-$solutiuon of the equation in 
\req{a17}
is uniquely defined. Namely we use that branch (of the three solutions) 
that is real for all positive values of $s$, and $s>0$ correponds with 
$w\in\RR$.

After these preparation we obtain the standard form (cf. \req{a1})
$$
U\left(-\brg12\mu^2,\mu t\sqrt{2}\right)e^{-zA}=
{\sqrt{2\pi}\,e^{\br12\mu^2t^2}}
\left(\br\mu 2\right)^{\br12\mu^2+\br12}F_\eta(z),\en{a20}$$
where
$$F_\eta(z)=\br1{2\pi i}\int_{{\cal
C}}e^{z(\br13w^3-\eta w)}f(w)\,dw,\sq f(w)=\br1{\sqrt{s}}\br{ds}{dw}.$$
Taking into account the mapping in \req{a17}, we have
$$\br{ds}{dw}=s\br{ w^2-b^2 }{s^2-2ts+1},\sq f(w)=\sqrt{s}\br{w^2-b^2}
{s^2-2ts+1},\sq \eta=b^2.\en{a21}$$

As explained in the previous section, for the computation of the coefficients
$\alpha_n, \beta_n$, we need the coefficients $p_k^{(1)}, p_k^{(2)}$ of 
the expansions   (cf. \req{a9})
$$f(w)=\sk p_k^{(1)}(w-b)^k,\sq f(-w)=\sk p_k^{(2)}(w-b)^k.\en{a22}$$
It turns out that $p_0^{(1)}=p_0^{(2)}$. Indeed, consider the expansions:
$$s=s_++\sum_{k=1}^\infty s_k^+(w-b)^k,\sq s=s_-+\sum_{k=1}^\infty 
s_k^-(w+b)^k.\en{a23}$$
Using the expression of $ds/dw$ in \req{a21} and l'Hospital rule we obtain
$$
s_1^+=s_+\br{2b}{2(s_+-t)\,s_1^+},\sq\sq {\rm so\ that}\sq
s_1^+=\br{\sqrt{bs_+}}{(t^2-1)^\br14}=\sqrt{\br{bs_+}{\sinh\theta}}.
$$
The square root has the plus sign because $ds/dw$ is positive if $w\in\RR$,
as follows from the first relation in \req{a21} and the properties of the 
mapping. From the expression \req{a20a} for $f(w)$ we obtain:
$$
p_0^{(1)}=f(b)=\br{s_1^+}{\sqrt{s_+}}=\sqrt{\br{b}{\sinh\theta}}.
$$
Analogously,
$$
s_1^-=\sqrt{\br{bs_-}{\sinh\theta}} \sq\sq {\rm and} \sq\sq
p_0^{(2)}=f(-b)=\sqrt{\br{b}{\sinh\theta}}=p_0^{(1)}. \en{p10p20}
$$

In order to avoid expressions with algebraic roots in the computations,
it is convenient to consider expansions like \req{a22} for the function
$\tilde{f}(w)=f(w)/f(b)$. The corresponding coefficients we denote by
$\tilde{p}_k^{(1)}$ and $\tilde{p}_k^{(2)}$. Besides, to avoid algebraic
roots in the expansions of \req{a23} we replace $t$ by a new variable
$$
u= \sqrt{2b} \left(\br{t-1}{t+1}\right)^{\br14}, \sq\sq
{\rm so\ that} \sq t=\br{4b^2+u^4}{4b^2-u^4}. 
$$
Then
$$
s_+=\br{2b+u^2}{2b-u^2}, \sq s_1^+=\br{2b+u^2}{2u}.
$$
Other coefficients $s_k^+$ can be obtained by deriving a reccurence relation
for them from the differential equation in \req{a21}. They are rational
functions in $u$ and $b$. The coefficients $s_k^-$ can be obtained 
from the corresonding $s_k^+$ by changing the sign of both $u$ and $b$.
In particular, 
$$
s_-=\br{2b-u^2}{2b+u^2}, \sq s_1^-=\br{2b-u^2}{2u}.
$$

Further, the coefficients $\tilde{p}_k^{(1)}$ and $\tilde{p}_k^{(2)}$
can be computed using
$$
f(w)=\br1{\sqrt{s}}\br{ds}{dw}=2\br{d\sqrt{s}}{dw}.
$$
Recall that $\sqrt{s}$ satisfies the differential equation
$2s\,dS/dw=S\,ds/dw$. It is convenient to compute the power series
(in $w-b$) solution $S_+$ of this equation with $S_+(b)=4u/(2b-u^2)$. Then
$\tilde{f}(w)=dS_+/dw$ and the coefficients $\tilde{p}_k^{(1)}$ are obtained
easily. The coefficients $\tilde{p}_k^{(2)}$ can be obtained by 
changing the sign of both $b$ and $u$ in the expression for
$(-1)^k\tilde{p}_k^{(1)}$.

Application of the algorithm of the  previous section gives the coefficients
$\alpha_j$, $\beta_j$ for the expansion of $\tilde{f}(w)$, and these
coefficients are rational functions in $b$ and $u$. We write them 
in a more compact form as rational functions in $\eta=b^2$ and
$$
\xi=\br{u^4+4b^2}{4u^2}=\br{bt}{\sqrt{t^2-1}}.
$$
Then first few coefficients in the expansion \req{a5} are:
$$\el{
\alpha_0=1, \sq\sq & \beta_0=0,\cr
\alpha_1=\br{1}{48}, \sq\sq  &
\beta_1=\br{5\xi^3-6\eta\xi-5}{48\,\eta^2},
}$$
$$\el{
\alpha_2=\br{385\xi^6-924\eta\xi^4+684\eta^2\xi^2-143\eta^3
+70\xi^3-84\eta\xi-455}{4608\,\eta^3}, \cr 
\beta_2=\br{\beta_1}{48},\sq\sq
\alpha_3=\br{\alpha_2}{48}-\br{2021}{34560}\alpha_1,\cr
\beta_3=\br{425425\xi^9-1531530\eta\xi^8+2040012\eta^2\xi^5
-28875\xi^6-1189005\eta^3\xi^3+69300\eta\xi^4}{3317760\,\eta^5}\cr
+\br{259110\eta^4\xi-51300\eta^2\xi^2+28875\xi^3
+10725\eta^3-34650\eta\xi-425425}{3317760\,\eta^5},\cr
\beta_4=\br{\beta_3}{48}-\br{2021}{34560}\beta_2.\cr
}
$$
The linear relations between the coefficients follow from
expansion (8.11) in \ref{Olver}(1959), where both power
series factors of $\Ai$ and $\Ai'$ contain only even powers of our $z$
(in Olver's notation, $z=\br12\mu^2$), but the whole expansion is
multiplied by function $g(z)$ with known asymptotics.
Olver also notes that coefficients in the asymptotic expansion
of $U(a,x)$ in terms of Airy functions can be lineraly determined
from the asymptotic expansion (of the same function) in terms of
elementary functions; see formulas (8.12), (8.13) in \ref{Olver}(1959).

The coefficients $\alpha_n, \beta_n$ are analytic functions at
$\eta=0$ and we can expand them in Maclaurin series. 
The first few coefficients are expanded as follows:
$$\el{
\beta_1=-\br{9}{560}+\br{7}{1800}\eta-\br{1359}{1078000}\eta^2
+\br{7}{16250}\eta^3-\br{152723}{1018710000}\eta^4+
\br{3997}{75968750}\eta^5+\ldots,\cr
\alpha_2=-\br{199}{115200}+\br{6849}{4928000}\eta
-\br{737}{1040000}\eta^2+\br{46711}{142560000}\eta^3
-\br{975823}{6806800000}\eta^4+\ldots.\cr
}
$$
The radius of convergence equals $(3\pi/2)^{2/3}=2.81...$. 
This number follows from the singularity of the mapping
given in (19), with $\theta$ defined in (18). 
The mapping is singular at $t=-1$.

\sect{Maple code}

General case. For an input one has to (re)define functions {\sffont AiryPw}
and {\sffont AiryPm} specifying the coefficients in \req{a9}.
Output is given by functions {\sffont AiryAlpha} and {\sffont AiryBeta},
which return the coefficients in \req{a5}. For convenience, one may rename
the global variable {\sffont AiryB} using {\sffont alias}.

\vskip 6pt
{\bfont \parindent=14pt \parskip=-1pt  
AiryAlpha:= proc( n )  normal(AiryGamma(n,0))  end:

AiryBeta:=  proc(n)  normal(AiryDelta(n,0))  end:

AiryGamma:= proc( n, k )
\maplenl if n=0 then  AiryC(k)
\maplenl else factor(
   (2*k+1)*AiryDelta(n--1,k+1)+2*AiryB\^{}2*(k+1)*AiryDelta(n--1,k+2) )
\maplenl fi
\par end:

AiryDelta:= proc(n, k )
\maplenl if n=0 then  AiryD(k)
\maplenl else 2*(k+1)*AiryGamma(n--1,k+2)
\maplenl fi
\par end:

AiryC:= proc( k )
local j;
\maplenl if k=0 then  AiryFe(0)
\maplenl else factor(
\maplenl\ \ sum( '(--1)\^{}(k--j)*j/(2*k--j)*binomial(2*k--j,k)/(2*AiryB)\^
{}(2*k--j)*AiryFe(j)', 'j'=1..k ) )
\maplenl fi
\par end:

AiryD:= proc( k )
local j;
\maplenl if k=0 then  AiryFoe(0)
\maplenl else factor(
\maplenl\ \  sum( '(--1)\^{}(k--j)*j/(2*k--j)*binomial(2*k--j,k)/(2*AiryB)\^
{}(2*k--j)*AiryFoe(j)', 'j'=1..k ) )
\maplenl fi
\par end:

AiryFoe:= proc( k )
\maplenl if k$<$0 then  0
\maplenl else expand( (AiryFo(k)--AiryFoe(k--1))/AiryB );
\maplenl  fi
\par end:

AiryFe:= proc(k)  (AiryPw(k)+AiryPm(k))/2  end:

AiryFo:= proc(k)  (AiryPw(k)--AiryPm(k))/2  end:
\vskip 4pt}

To find the coefficients of the expansion of parabolic cylinder function
$U(a,x)$ (up to the multiple in \req{p10p20}) one has to assign
{\par \parindent=14pt
\bfont AiryPw:= ParCyPw; \ \ AiryPm:= ParCyPm;
\par \noindent}The global variables are {\sffont AiryB, ParCyU, ParCyXi},
they correspond to variables $b,u,\xi$ in the text. The coefficients in $b$
and $u$ would be returned by {\sffont AiryAlpha} and {\sffont AiryBeta},
and coefficients in $b$ and $\xi$ --- by {\sffont ParCyAlpha} and
{\sffont ParCyBeta}.

\vskip 6pt
{\bfont \parindent=14pt \parskip=-1pt
alias( ParCyUa=RootOf( z\^{}4--4*ParCyXi*z\^{}2+4*AiryB\^{}2, z)):
\par\# Algebraic relation between ParCyU and ParCyXi

ParCyAlpha:=  proc(k)  factor( evala(subs(ParCyU=ParCyUa,AiryAlpha(k))) ) end:

ParCyBeta:= proc(k)  factor( evala(subs(ParCyU=ParCyUa,AiryBeta(k))) ) end:

ParCySw:= proc( k )
option remember;
local T, s, a, w;
\maplenl if k=0 then  (2*AiryB+ParCyU\^{}2)/(2*AiryB--ParCyU\^{}2)
\maplenl elif k=1 then  (2*AiryB+ParCyU\^{}2)/2/ParCyU
\maplenl else
  s:= sum('a[i]*w\^{}i', 'i'=0..k);
\maplenl\ \ T:= coeff( expand(
(s\^{}2+2*(ParCyU\^{}4+4*AiryB\^{}2)/(ParCyU\^{}4--4*AiryB\^{}2)*s+1)
\maplenl\ \ \ \ *diff(s,w)--w*(w+2*AiryB)*s ), w, k);
\maplenl\ \ sort( factor( solve(
    subs( seq(a[i]=ParCySw(i),i=0..k--1), T), a[k])), ParCyU)
\maplenl fi
\par end:

ParCySm:= proc(k)
  subs( ParCyU=--ParCyU, AiryB=--AiryB, ParCySw(k) )
end:

ParCySqrtS:= proc( k )
option remember;
local j;
\maplenl if k=0 then  4*ParCyU/(2*AiryB--ParCyU\^{}2)
\maplenl else
    factor(  sum('(3/2*j--k)*ParCySw(j)*ParCySqrtS(k--j)','j'=1..k)
        /ParCySw(0)/k )
\maplenl fi
\par end:

ParCyPw:= proc(k)  (k+1)*ParCySqrtS(k+1)  end:

ParCyPm:= proc(k)
  (--1)\^{}k*subs(AiryB=--AiryB, ParCyU=--ParCyU, ParCyPw(k))
end:
}

\sect{References}
\def\ref#1{{\reffont#1\ }\advance \refnonum by 1}
\def\refno{{	{\rm \the\refnonum}}}
\newcount\refnonum
\refnonum=1
\parindent=18pt
\frenchspacing  \baselineskip=12pt \parskip 4pt plus 1pt

\item{[\refno]}
\ref{M. Abramowitz and I.A. Stegun} (1964), 
{\sl Handbook of mathematical functions with formulas, graphs and mathematical
tables}, Nat. Bur. Standards Appl. Series, {\bf 55}, 
U.S. Government Printing Office, Washington, D.C. (paperback edition 
published by Dover, New York). 

\item{[\refno]}
\ref{Bleistein, N.}
(1966), 
\rm{Uniform asymptotic expansions of integrals 
with stationary points and algebraic singularity}, 
\sl Comm. Pure Appl. Math., 
\bf 19,
\rm 353--370. 

\item{[\refno]}
\ref{Chester, C., Friedman, B. \& Ursell, F.}
(1957), 
\sl An extension of the method of steepest descent, 
\rm Proc. Cambridge Philos. Soc., 
\bf 53
\rm 599--611.

\item{[\refno]}
\ref{Dunster, T.M.}
(1989), 
\rm{Uniform asymptotic expansions for Whittaker's confluent hypergeometric
functions},    
\sl SIAM J. Math. Anal.,
\bf 20,
\rm 744--760.

\item{[\refno]}
\ref{Olde Daalhuis, A.B. \& Temme, N.M.}
(1994),
{\rm Uniform Airy type expansions of integrals},
{\sl SIAM J. Math. Anal.},
\bf 25,
\rm 304--321.

\item{[\refno]}
\ref{Olver, F.W.J.} (1959),
Uniform asymptotic expansions for
Weber para\-bolic cylinder functions of large orders,
{\sl J. Research NBS}, {\bf 63B}, 131--169.

\item{[\refno]}
\ref{Olver, F.W.J.} 
(1974),
{\sl Asymptotics and Special Functions}, 
Academic Press, New York. Reprinted in 1997 by A.K. Peters. 

\item{[\refno]}
\ref{Temme, N.M.}
\rm (1996). 
\sl Special functions: An introduction to the classical
functions of mathematical physics.
\rm Wiley, New York.

\item{[\refno]}
\ref{Temme, N.M.}
\rm (2000).
\rm{Numerical and asymptotic aspects of parabolic cylinder functions},
{\sl J. Comp. Appl. Math.},
\bf 121,
\rm 221--246.

\item{[\refno]}
\ref{Wong, R.} 
\rm (1989),
\sl Asymptotic approximations of integrals,
\rm Academic Press, New York.

\bye